\documentclass[twocolumn]{autart}

\usepackage{amsmath}
\usepackage{amsfonts}

\newcommand{\cL}{\mathcal{L}}
\newcommand*{\Set}[1]{\mathbb{#1}}
\newcommand{\sC}{\Set{C}}

\newcommand{\sN}{\Set{N}}
\newcommand{\sR}{\Set{R}}
\newcommand{\bg}{\bar{g}}
\newcommand{\bd}{\bar{d}}
\newcommand{\ceil}[1]{\left\lceil #1 \right\rceil}
\newcommand{\inprd}[1]{\left\langle #1 \right\rangle}

\DeclareMathOperator{\T}{T}
\DeclareMathOperator{\tr}{tr}
\DeclareMathOperator{\IQ}{IQ}
\DeclareMathOperator{\Ideal}{Ideal}
\DeclareMathOperator{\Qmod}{Qmod}

\DeclareMathOperator{\rank}{rank}

\DeclareMathOperator{\coef}{coef}
\DeclareMathOperator{\st}{s.t.}
\DeclareMathOperator{\diff}{d}

\begin{document}

\begin{frontmatter}
\title{Existence of Solutions for the Continuous Algebraic Riccati Equation 
via Polynomial Optimization \thanksref{footnoteinfo}}

\thanks[footnoteinfo]{This paper was not presented at any IFAC meeting. 
Corresponding author Juan Zhang.}

\author[zj]{Juan Zhang}\ead{zhangjuan@xtu.edu.cn},
\author[zwj]{Wenjie Zhao}\ead{zhaowj@hnfnu.edu.cn}

\address[zj]{Key Laboratory of Intelligent Computing and 
Information Processing of Ministry of Education, 
Hunan Key Laboratory for Computation and Simulation in Science and Engineering, 
School of Mathematics and Computational Science, 
Xiangtan University, Xiangtan 411105, Hunan, China}   
\address[zwj]{School of Mathematics and Statistics, 
Hunan First Normal University, Changsha 410205, Hunan, China \& 
School of Mathematics and Computational Science, 
Xiangtan University, Xiangtan 411105, Hunan, China}


\begin{abstract}
This paper studies the solution existence of the continuous-time algebraic
Riccati equation (CARE). We formulate the CARE as two constrained polynomial 
optimization problems, and then use Lasserre's hierarchy of semi-definite 
relaxations to solve them. Compared to the existing work, our approaches 
can obtain the exact positive semi-definite solution of the CARE 
when the coefficient matrices are symmetric. Moreover, 
our methods can detect the nonexistence of the positive semi-definite 
solution for the CARE. Numerical examples show that our approaches 
are effective compared to the existing methods. 
\end{abstract}

\end{frontmatter}

\section{Introduction}
Consider the continuous algebraic
Riccati equation (CARE) 
\begin{align} \label{CARE:general}
  A^{\T}P + PA - PRP + Q = 0,
\end{align}
where $R, Q \in \sR^{n \times n}$ are symmetric. 
The symbol $^{\T}$ denotes the transpose operator. One can refer to 
\cite{Lancaster1995ARE} for more details.

The CARE can be used in various areas, such as 
control system (\cite{Bini2012Book} and \cite{Hesp2018Book}), 
differential games \cite{Delf2007DiffGames}. 
Consider the linear quadratic optimal control system defined by
\begin{equation} \label{LS}
\begin{cases}
   \dot{x}(t) = Ax(t) + Bu(t), \quad x(0) = x_{0}, \\
   y(t) = Cx(t), 
\end{cases}
\end{equation}
where $x(t) \in \sR^{n}$, $u(t) \in \sR^{m}$, $y(t) \in \sR^l$ are the state vector, 
control (input) vector and output vector, respectively. $B \in \sR^{n\times m}$ 
and $C \in \sR^{l \times n}$. $x_0$ is a given initial state vector, and
\[
  \dot{x}(t) \,:=\, \frac{\diff{x}}{\diff{t}} 
  \,:=\, 
  \big( x_1'(t), ~x_2'(t),\dots,x_n'(t)\big)^{\T}.
\] 
The linear state feedback control law of system \eqref{LS} is
\[
u(t) = -Lx(t),
\]
with $L\in\sR^{m\times n}$ is the state feedback gain matrix.
Therefore, the closed-loop system has the following form
\[
  \dot{x}(t) = (A -BL)x(t), \quad x(0) = x_0.
\]
The optimal control means finding a control vector $u(t)$ among the allowed controls 
to minimize the following performance index
\[
  J = \int_{0}^{\infty}\big(x(t)^{\T}Qx(t) + u(t)^{\T}Wu(t)\big)\,\mathrm{d}t,
\]
where $Q \in \sR^{n \times n}$, $W \in \sR^{m \times m}$ are weighting matrices.
$Q=C^{\T}C$ is positive semi-definite and $W$ is positive definite.
The optimal feedback gain matrix is expressed by
\[
  L^* = W^{-1}B^{\T}P,
\]
where $P$ is the solution of the CARE \eqref{CARE:general}, 
$R = BW^{-1}B^{\T}=R_1R_1^{\T}$ is positive semi-definite.
Refer to \cite{Hesp2018Book} for more details.

In general, there are many solutions to the CARE. One is interested 
in the stabilizing solution,  which means that the state vector $x(t)$ 
will remain bounded or converge to a steady state. If $P$ is a stabilizing solution, 
then the closed-loop system $\dot{x}(t)=(A-RP)x(t)$ is stable 
(i.e., all eigenvalues of $A-RP$ have negative real parts).

There exist some works involving solution existence conditions 
for the CARE when $R$ and $Q$ are positive semi-definite 
(\cite{Kucera1972ExistenceA} and \cite{Kano1987Existence}).
Among these, sufficient and necessary conditions are provided 
for the existence of the stabilizing solutions or the (semi)definite solutions, 
given in terms of controllability (or stabilizability) and observability 
(or detectability) of control systems. Note that stabilizability and 
detectability are weaker than controllability and observability, 
respectively. Thus, detecting the solution existence of the
CARE via other methods is a challenge. 

Assume that the CARE has a positive semi-definite solution, many numerical 
methods have been developed. The standard computational methods include 
Schur method \cite{Laub1979Schur}, structure-preserving doubling algorithm 
(\cite{Chu2005StruPresDoubling} and \cite{Lin2006StruPresDoubling}), 
matrix sign function method \cite{Byers1987MatrixSignFunction}, 
Newton-based method (\cite{Benner1998Newton} and \cite{Bini2008Newton}), 
inexact Newton method \cite{Wang2014InexactNewton}, etc. 
The direct applications of low-rank subspace alternating direction implicit iteration 
have also been recently studied (\cite{Benner2013LowRankADI} and \cite{Benner2018LowRankADI}). 
Another approach is based on projection onto Krylov-type subspace method using 
standard, rational and extended block Arnoldi processes (\cite{Jbilou2003Krylov} and \cite{Simon2016RationalKrylow}).

Viewing these numerical methods, Schur method does not fully use
the Hamiltonian matrix characteristics. Krylov subspace method is inaccurate 
as it uses Galerkin orthogonality conditions to obtain an approximate solution 
at every low-order Krylov subspace. The Newton method has an extensive computation. 
Therefore, designing efficient and stable algorithms for the CARE is another challenge. 

In this paper, we focus on computing the positive semi-definite 
solution of the CARE by solving polynomial optimization problems. 
Our contributions are summarized as follows.
\begin{itemize}

\item
We formulate the CARE as two polynomial optimization problems, 
one exploiting Schur decomposition of a positive semi-definite matrix 
involves standard polynomial constraints, and the other incorporates 
polynomial matrix inequality constraint. 
We can obtain a positive semide-finite solution of the CARE 
when $R$ and $Q$ are symmetric by solving polynomial optimization problems.
 
\item
In general, Lasserre's hierarchy of semi-definite relaxations 
is used to obtain the exact positive semi-definite solution of 
the CARE with a small relaxation order 
by solving polynomial optimization problems. 
Moreover, if the positive semi-definite solution has an upper bound, 
our methods are convergent without the assumption of the archimedean condition.

\item
Moreover, we can detect the nonexistence of the positive semi-definite 
solution of the CARE by solving polynomial optimization 
problems instead of checking the existence conditions of 
the positive semi-definite solution. 

\end{itemize}

The rest of the paper is organized as follows.
In Section 2, we recall some basic theories of polynomial optimization.
In Section 3, we formulate the CARE as a constrained polynomial 
optimization problem with equality and inequality constraints, 
and Lasserre's hierarchy of semi-definite relaxations 
is used to solve constrained polynomial optimization problem 
for obtaining a positive semi-definite solution of the CARE.
In Section 4, we present another algorithm for getting 
the positive semi-definite solution of the CARE by 
solving a polynomial matrix optimization problem.
In Section 5, numerical examples are given to illustrate that 
our algorithms are effective. Finally, we draw some concluding 
remarks and prospects in Section 6.

\section{Preliminaries}

\subsection{Notation}

The symbol $\sN$ (resp., $\sR$, $\sC$) denotes the set of nonnegative integral
(resp., real, complex) numbers.  Denote $\sR[x] := \sR[x_1, \cdots , x_l]$
by the ring of polynomials in $x := (x_1,\cdots, x_l)$ with real coefficients.
The $\sR[x]_d$ stands for the set of polynomials in $\sR[x]$ with degrees $\leq d$.
For a polynomial $p$, $\deg(p)$ denotes its degree.
For $t\in \sR$, $\lceil t \rceil$ denotes the smallest integer $\geq t$.
$I_n$ is an $n\times n$ identity matrix.
For $\alpha := (\alpha_1, \cdots, \alpha_l) \in \sN_l$ with an integer $l > 0$,
$|\alpha| := \alpha_1 + \alpha_2 + \cdots + \alpha_l$. 
For an integer $d > 0$, denote
\[
\sN^l_d:= \{\alpha \in \sN^l~|~|\alpha| \leq d\}.
\]
For $z=(x_1,\cdots,x_l)$ and $\alpha=(\alpha_1,\cdots,\alpha_l)$, denote
\begin{gather*}
x^\alpha:= x^{\alpha_1}_1 \cdots x^{\alpha_l}_l, \\
[x]_d := [1 ~~ x_1 ~~\cdots ~~ x_n ~~ x^2_1 ~~ x_1x_2 ~~\cdots~~ x^d_l]^{\T}.
\end{gather*}
$X \succeq 0$ means that $X$ is a symmetric positive semi-definite matrix.
For a vector $a$, $\| a \|$ denotes standard Euclidean norm of $a$.
For a matrix $A$, $\| A \|_F$ means Frobenius norm of $A$.

\subsection{Ideals, sum of squares, quadratic modules, localizing and moment matrices}

In this subsection, we review some basic results in the polynomial ring $\sR[x]$.
An ideal $I$ of $\sR[x]$ is a subset such that $I \cdot\sR[x] \subseteq I$
and $I +I\subseteq I$. For a tuple $p = (p_1, \cdots, p_k)$ of polynomials in $\sR[x]$,
Ideal($p$) denotes the smallest ideal containing all $p_i$, which is the set
\[
\Ideal(p) := p_1\cdot\sR[x] + \cdots + p_k\cdot\sR[x].
\]
In computation, we often need to work with its truncation
\[
\Ideal(p)_{2k} := p_1\cdot \sR[x]_{2k-\deg(p_1)} + \cdots +
p_k\cdot \sR[x]_{2k-\deg(p_k)}.
\]

A polynomial $\sigma$ is said to be a sum of squares (SOS) 
if $\sigma = s^2_1+s^2_2+\cdots+s^2_k$ for some polynomials 
$s_1, s_2, \cdots, s_k$. Whether or not a polynomial is SOS 
can be checked by solving a semi-definite program (SDP) \cite{Lasserre2010Book}. 
Obviously, if a polynomial is SOS, then it is nonnegative. 
However, the reverse may not be true.
The set of all SOS polynomials in $\sR[x]$ is denoted as $\Sigma[x]$, 
and its $d$-th truncation is
\[
  \Sigma[x]_d := \Sigma[x] \cap \sR[x]_d.
\]
For a tuple $q = (q_1, q_2, \cdots, q_t)$ of polynomials in $\sR[x]$, 
its quadratic module is
\[
  \Qmod(q) := \Sigma[x] + 
  q_1 \cdot \Sigma[x] + 
  q_2 \cdot \Sigma[x] + \cdots + 
  q_t \cdot \Sigma[x].
\]
We often need to work with the truncation
\begin{equation*}
\begin{split}
  \Qmod(q)_{2k} := \Sigma[x]_{2k} +
  q_1 \cdot \Sigma[x]_{2k - \deg(g_1)} + \\
  q_2 \cdot \Sigma[x]_{2k - \deg(g_2)} + \cdots +
  q_t \cdot \Sigma[x]_{2k - \deg(q_t)}.
\end{split}
\end{equation*}

Let $x := (x_1, x_2, \cdots, x_n)$ be a vector. Denote $\sR^{\sN^n_d}$
by the space of real sequences indexed by $\alpha\in \sN^n_d$.
The vector $y := (y_\alpha)_{\alpha \in \sN^n_d}\in \sR^{\sN^n_d}$
is called a truncated multi-sequence (tms) of degree $d$.
It gives the Riesz functional $\cL_y$ acting on $\sR[x]_d$ as
\[
  \cL_y \left( \sum_{\alpha\in \sN^n_d}f_\alpha x^\alpha \right)
  :=\sum_{\alpha \in \sN^n_d}f_\alpha y_\alpha.
\]
For $f\in \sR[x]_d$ and $y\in \sR^{\sN^n_d} $, we denote
\[
  \inprd{f, y} := \cL_y(f).
\]
Consider a polynomial $q \in \sR[x]_{2k}$ with $\deg(q) \leq 2k$.
The $k$-th localizing matrix of $q$, generated by a tms $y \in \sR^{\sN^n_{2k}}$
is the symmetric matrix $L_q^{(k)}(y)$ such that
\[
\coef(a_1)^T L_q^{(k)}(y) \coef(a_2) = \cL_y(qa_1a_2)
\]
for all $a_1, a_2 \in \sR[x]_{k-\lceil \deg(q)/2\rceil}$,
in which $\coef(a_i)$ denotes the coefficient vector of $a_i$.

When $q = 1$ (the constant one polynomial), $L_1^{(k)}(y)$ is called a moment matrix.
We use $M_k(y):=L_1^{(k)}(y)$. The columns and rows of $L_q^{(k)}(y)$,
as well as $M_k(y)$, are labelled by $\alpha\in \sN^n$ with
$2|\alpha| \leq 2k-\deg(q)$.

Let $H$ be an $l\times l$ symmetric matrix, whose each element
$H_{ij}$ is a polynomial in $\sR[x]$. The $k$-th localizing
matrix of $H$ generated by the tms $y \in\sR^{\sN_{2k}^n}$ is
the block symmetric matrix $L_{H}^{(k)}(y)$, which is defined as
\[
L_{H}^{(k)}(y) := \big(L_{H_{ij}}^{(k)}(y)\big)_
{1 \leq i \leq l, 1 \leq j \leq l},
\]
and each block $L_{H_{ij}}^{(k)}(y)$ is a standard localizing matrix
of the polynomial $H_{ij}$.

Let $\Sigma[x]^{l\times l}$ be the cone of all sums of
$s_1s_1^{\T}+s_2s_2^{\T}+\cdots+s_rs_r^{\T}$,
where $s_1,s_2,\cdots,s_r\in\sR[x]^l$.
When $l=1$, $\Sigma[x]^{l\times l}$ is $\Sigma[x]$.
The quadratic module of $H$ is
\[
\Qmod(H) := \Sigma[x] + \left\{ \tr(HS):
S \in \Sigma[x]^{l\times l} \right\}.
\]
The $k$-th truncation of $\Qmod(H)$ is defined as
\begin{multline*}
\Qmod(H)_{2k} := \Sigma[x]_{2k} + \\ 
\left\{ \tr(HS) ~\left|~
\begin{array}{l}
	S \in \Sigma[x]^{l \times l}, \\
	\deg(H_{ij}S_{ij}) \leq 2k, \\
	\forall~1 \leq i \leq l,1 \leq j \leq l
\end{array} \right. \right\}.
\end{multline*}

For two tuples $p = (p_1, \cdots, p_k)$, $q = (q_1, \cdots, q_t)$
of polynomials in $\sR[x]$, and a symmetric matrix $H$, we denote
\begin{align*}
  \IQ(p, q, H) :=& \Ideal(p) + \Qmod(q)+ \Qmod(H),\\
  \IQ(p, q, H)_{2k} :=& \Ideal(p)_{2k} + \Qmod(q)_{2k}+ \Qmod(H)_{2k}.
\end{align*}
The set $\IQ(p, q, H)$ (resp., $\IQ(p, q, H)_{2k}$) is a convex cone
that is contained in $\sR[x]$ (resp., $\sR[x]_{2k}$).

The set $\IQ(p, q, H)$ is said to be archimedean if there exists
$\sigma \in \IQ(p, q, H)$ such that $\sigma(x) \geq 0$.
If $\IQ(p, q, H)$ is archimedean, then the set
$K := \{ p(x) = 0, q(x) \geq 0, H \succeq 0 \}$
must be compact. The reverse is not always true.
However, if $K$ is compact, say, $K \subseteq B(0, \tilde{R})$
(the ball centered at $0$ with radius $\tilde{R}$), then 
$\IQ(p, \tilde{q}, H)$ with $\tilde{q} = (q, \tilde{R}^2 - \|x\|^2)$ 
is always archimedean.
Under the assumption that $\IQ(p, q, H)$ is archimedean,
every polynomial in $\sR[x]$, which is strictly positive on $K$, 
must belong to $\IQ(p, q, H)$, which is the so-called 
Putinar’s Positivstellensatz (\cite{Putinar1993Qmod} and \cite{Nie2023Book}).
Interestingly, under some optimality conditions, 
if a polynomial is nonnegative (but not strictly positive) over $K$, 
then it belongs to $\IQ(p, q, H)$. Please see \cite{Nie2023Book} 
and \cite{Nie2015LocalMinimums} for more details.

\section{A semi-definite relaxation algorithm for the CARE}

In general, the CARE \eqref{CARE:general} may have several solutions,
such as negative (semi)definite solutions, indefinite solutions 
or positive (semi)definite solutions. For instance, let 
\[
  A = \begin{pmatrix}  -1 & 0 \\ 0 & -2  \end{pmatrix}, \quad
  R = \begin{pmatrix}   1 & 0 \\ 0 &  2  \end{pmatrix}, \quad
  Q = \begin{pmatrix}   2 & 0 \\ 0 &  3  \end{pmatrix} 
\]
in the CARE. It has one negative definite solution, 
two indefinite solutions and one positive definite solution below
\begin{gather*}
  \begin{pmatrix}  -2.7321 & 0 \\ 0 & -2.5811  \end{pmatrix}, \quad
  \begin{pmatrix}  -2.7321 & 0 \\ 0 &  0.5811  \end{pmatrix}, \\
  \begin{pmatrix}   0.7321 & 0 \\ 0 & -2.5811  \end{pmatrix}, \quad
  \begin{pmatrix}   0.7321 & 0 \\ 0 &  0.5811  \end{pmatrix}.
\end{gather*}
In this paper, we focus on finding a positive semi-definite solution. 
Therefore, an additional positive semi-definite 
constraint for matrix $P$ is added. Since any real positive semi-definite 
matrix has a Cholesky decomposition, then we have
\begin{equation} \label{CD}
	P = XX^{\T}
\end{equation}
where 
\begin{align*}
	X = \begin{pmatrix}
X_{11} & 0      & 0      & \cdots & 0 \\
X_{21} & X_{22} & 0      & \cdots & 0 \\
\vdots & \vdots & \vdots & \ddots & 0 \\
X_{n1} & X_{n2} & 0      & \cdots & X_{nn} \\
	\end{pmatrix}\in \sR^{n\times n}
\end{align*}
is a lower triangular matrix with nonnegative diagonal elements (Cholesky factor). 
Applying Cholesky decomposition \eqref{CD} to the CARE \eqref{CARE:general}, we get 
\begin{equation} \label{CARE:CD}
	A^{\T} XX^{\T} + XX^{\T} A - XX^{\T} R XX^{\T} + Q = 0.
\end{equation}
For convenience, denote
\begin{align*}
	& f(X) := \sum_{i,j=1}^{n}b_{ij}X^2_{ij}+
	\sum_{\substack{i,j,k,l=1, \\ i\neq k,j\neq l}}^{n}c_{ij}X_{ij}X_{kl}+
	\sum_{i,j=1}^{n}d_{ij}X_{ij},  \\
	& g_{ij}(X) := (A^{\T}XX^{\T}+XX^{\T}A-XX^{\T}RXX^{\T}+Q)_{ij}, \\
    & \qquad 1 \leq j \leq i \leq n, \\
	& g := (g_{ij})_{1 \leq j \leq i \leq n}, \quad
	h := (X_{11}, X_{22}, \dots, X_{nn}),
\end{align*}
where all $b_{ij},c_{ij},d_{ij}$ are random positive constants.

The polynomial system of \eqref{CARE:CD} can be described
as the following optimization problem
\begin{align}\label{opt:CARE:CD}
	\begin{cases}
\min & f(X),\\
\st  & g(X) = 0, \\
& h(X) \geq 0,
	\end{cases}
\end{align}
which is a standard optimization problem with equality  
and inequality constraints.

Let
\[
d_0 := \ceil{\max \{ \deg(g), ~\deg(h) \}/2}=2.
\]
Apply Lasserre’s hierarchy of semi-definite relaxations
to the problem \eqref{opt:CARE:CD}, for integers $k \geq \max(1,d_0)=2$,
we get the following $k$-th SDP problem 
($k$ is called as relaxation order)
\begin{align}\label{opt:CARE:CD:Relax}
	\begin{cases}
\inf_y & \inprd{f, y}, \\
\st    & M_k(y) \succeq 0, \quad y_0 = 1, \\
& L_{g_{ij}}^{(k)}(y) = 0, ~1 \leq j \leq i \leq n, \\
& L_{h_{j}}^{(k)}(y)\succeq 0, ~j=1, 2, \cdots, n. \\
	\end{cases}
\end{align}
Let $f^*$ and $f_k$ be the optimal value of the problem \eqref{opt:CARE:CD}
and the relaxation \eqref{opt:CARE:CD:Relax}, respectively.
Because the problem \eqref{opt:CARE:CD:Relax} is a relaxation problem 
of the problem \eqref{opt:CARE:CD}, 
we have $f_k \leq f^*$ for each relaxation order $k$.
As the relaxation order $k$ increases, the feasible set 
of the problem \eqref{opt:CARE:CD:Relax} will shrink. 
Hence, we have the monotonically increasing sequence
\[
  f_k ~\leq~ f_{k+1} ~\leq~ f_{k+2} ~\leq~ \cdots ~\leq~ f^*.
\]
Suppose $y^*$ is an optimal solution of the problem \eqref{opt:CARE:CD:Relax}. 
If the flat truncation condition holds \cite{Nie2013FlatTruncation} for $y^*$, 
that is, there exists an integer $t$ and $2 \leq t \leq k$ such that
\begin{align} \label{cdnt:cd:ft}
  \rank \big( M_t(y^*) \big) = \rank \big( M_{t-2}(y^*) \big),
\end{align}
then we can extract at least $\rank(M_t(y^*))$ minimizer 
of the problem \eqref{opt:CARE:CD} by the algorithm 
in \cite{Henrion2005extractSolution}. 
The flat truncation condition is an important property for checking 
the convergence of the relaxation problem \eqref{opt:CARE:CD:Relax}.
We now give the following algorithm for solving the problem \eqref{opt:CARE:CD}.
\begin{alg}
\label{algo:CARE:CD}

Let $f$, $g$ and $h$ be as in \eqref{opt:CARE:CD}
and $k = 2$, do the following:

\begin{itemize}

\item[1:]
Solve the SDP problem \eqref{opt:CARE:CD:Relax}.

\item[2:]
If the relaxation \eqref{opt:CARE:CD:Relax} is infeasible,
the CARE \eqref{CARE:general} has no positive semi-definite solution and stop.
Otherwise, solve it for a minimizer $y^*$. 
Let $t := 2$ and go to Step 3.

\item[3:]
Check whether or not $y^*$ satisfies the rank condition
\eqref{cdnt:cd:ft}. If it is satisfied, extract $\rank(M_t(y^*))$ minimizers 
for \eqref{opt:CARE:CD} and stop. Otherwise, go to Step 4.

\item[4:]
If $t < k$, let $t := t + 1$ and go to Step 3.
Otherwise, let $k := k + 1$ and go to Step 1.

\end{itemize}
\end{alg}

In Step 3, the flat truncation is a sufficient condition for obtaining 
the optimal solution of the problem \eqref{opt:CARE:CD}. 
Moreover, Algorithm \ref{algo:CARE:CD} is implemented by MATLAB 
software package {\tt GloptiPoly 3} \cite{Heniron2009GloptiPoly}.

The dual problem of the problem \eqref{opt:CARE:CD:Relax} is 
\begin{align}\label{opt:CARE:CD:Relax:Dual}
	\begin{cases}
\max & \gamma, \\
\st  & f-\gamma \in \IQ(g,h)_{2k},
	\end{cases}
\end{align}
which is also an SDP problem. Let $f_k^*$ be the optimal value of 
the problem \eqref{opt:CARE:CD:Relax:Dual}. According to weak duality, 
we have $f_k^* \leq f_k$ for each relaxation order $k$.
In addition, the feasible set of the problem \eqref{opt:CARE:CD:Relax:Dual}
will expand as the relaxation order $k$ increases. Therefore, 
we have the following monotonically increasing sequence
\[ 
  f_k^* ~\leq~ f_{k+1}^* ~\leq~ f_{k+2}^* ~\leq~ \cdots ~\leq~ f^*. 
\]
If the archimedean condition holds, then $f_k^*\to f^*$ as $k\to +\infty$. 
Please see \cite{Lasserre2001PolyOpt} for more details. When there is an interior point 
in the feasible set of the problem \eqref{opt:CARE:CD:Relax}, i.e., 
the  strong duality holds, then the problem \eqref{opt:CARE:CD:Relax:Dual}
achieve its optimal value and $f_k = f_k^*$ for each relaxation order $k$. 
Furthermore, If $f_k^*=f^*$ for some $k < +\infty$, we say the finite 
convergence holds for the problem \eqref{opt:CARE:CD:Relax:Dual}.

If the CARE \eqref{CARE:general} has an upper bound $P_u$ of a positive semi-definite 
solution $P_*$ under some conditions, then $P_u-P_*$ is positive semi-definite.
If the following constraint
\begin{align*} \label{care:X:upperbound}
  \| P_u \|_{F}^2 - \|XX^{\T}\|_{F}^2 \geq 0
\end{align*}
is added to the problem \eqref{opt:CARE:CD}, 
then the archimedean condition always holds.

For the problems \eqref{opt:CARE:CD} and \eqref{opt:CARE:CD:Relax}, 
we give the following results.
\begin{thm}
Let $f^*$ and $f_k$ be the optimal value of the problem \eqref{opt:CARE:CD}
and the relaxation \eqref{opt:CARE:CD:Relax}, respectively.
Assume that the archimedean condition holds, then
\begin{enumerate}
\item [\rm i)]
The relaxation \eqref{opt:CARE:CD:Relax} is infeasible for some $k$
if and only if the CARE \eqref{CARE:general} has no positive semi-definite solution.

\item [\rm ii)] If the relaxation \eqref{opt:CARE:CD:Relax} 
has a solution $y^*$ and the rank condition \eqref{cdnt:cd:ft} holds, 
then $f_k=f^*$ and there exists a positive semi-definite solution 
$X_*X_*^{\T}$ for the CARE \eqref{CARE:general} where $X_*$ is the optimal 
solution of the problem \eqref{opt:CARE:CD}.
\end{enumerate}
\end{thm}
	
\begin{pf}
i) Since the problem \eqref{opt:CARE:CD:Relax} is a relaxation problem of \eqref{opt:CARE:CD},
if the relaxation \eqref{opt:CARE:CD:Relax} is infeasible for some $k$, then the problem 
\eqref{opt:CARE:CD} is infeasible. Therefore, the CARE \eqref{CARE:general} has no positive 
semi-definite solution. 

Conversely, if the CARE \eqref{CARE:general} has no positive semi-definite solution,
then the set $\{ X\in\sR^{n\times n}~|~g(X) = 0, ~ h(X) \geq 0 \}$ is empty.
Hence, the problem \eqref{opt:CARE:CD} is infeasible.  
Because the archimedean holds, we have $-1 \in \IQ(g,h)$.
So $-1 \in \IQ(g, h)_{2k}$ for all such $k$ big enough, which means that
\eqref{opt:CARE:CD:Relax:Dual} is unbounded from above for all big $k$.
By weak duality, the relaxation \eqref{opt:CARE:CD:Relax} is infeasible.

ii) The result is from the classical achievements in
\cite{Nie2013FlatTruncation} and \cite{Curto2005TKM}.
\end{pf}

\section{Another polynomial optimization algorithm for the CARE}
In this section, another polynomial optimization algorithm with 
the positive semi-definite constraint is introduced 
for the CARE \eqref{CARE:general}. 

Define the following notation
\begin{align*}
	& \bar{f}(P):=\sum_{i,j=1}^{n}\bar{b}_{ij}P^2_{ij}+
	\sum_{\substack{i,j,k,l=1, \\ i \neq k,j \neq l}}^{n}\bar{c}_{ij}P_{ij}P_{kl}+
	\sum_{i,j=1}^{n}\bar{d}_{ij}P_{ij}, \\
	& \bg_{ij}(P):=\left(A^{\T}P+PA-PRP+Q\right)_{ij}, ~1 \leq i \leq j \leq n, \\
	& \bar{g} := (\bg_{ij})_{1 \leq i \leq j \leq n}, \quad
\end{align*}
where all $\bar{b}_{ij},\bar{c}_{ij},\bar{d}_{ij}$ are random positive constants.
We have the following polynomial optimization problem 
\begin{equation}\label{opt:CARE:PSD}
\begin{cases}
  \min & \bar{f}(P),  \\
  \st  & \bg(P) = 0,  \\
       & P \succeq 0.
\end{cases}
\end{equation}
Let
\[
  \bd_0 := \ceil{\max \{ \deg(\bg),~\deg(P) \}/2}=1.
\]
We apply Lasserre’s hierarchy of semi-definite relaxations
to solve the problem \eqref{opt:CARE:PSD}. 
For integers $k \geq \max(1,\bd_0)=1$,
the $k$-th semi-definite relaxation is
\begin{align}\label{opt:CARE:PSD:Relax}
	\begin{cases}
\inf_{\bar{y}} & \inprd{\bar{f}, \bar{y}}, \\
\st & M_k(\bar{y})\succeq 0, ~y_0=1, \\
& L_{\bg_{ij}}^{(k)}(\bar{y}) = 0, ~1 \leq i \leq j \leq n, \\
& L_P^{(k)}(\bar{y}) \succeq 0. 
	\end{cases}
\end{align}
Let $\bar{f}^*$ and $\bar{f}_k$ be the optimal value of the problem
\eqref{opt:CARE:PSD} and \eqref{opt:CARE:PSD:Relax}, respectively.
Since the problem \eqref{opt:CARE:PSD:Relax} is a relaxation problem 
of the problem \eqref{opt:CARE:PSD} with relaxation order $k$, we have 
$\bar{f}_k \leq \bar{f}^*$ for each $k$. 
Moreover, the following monotonically increasing sequence 
\[
 \bar{f}_k ~\leq~ \bar{f}_{k+1} ~\leq~ \bar{f}_{k+2} ~\leq~ \cdots ~\leq~ \bar{f}^*.
\]
Assume that $\bar{y}^*$ is a minimizer of the problem \eqref{opt:CARE:PSD:Relax}. 
If $\bar{y}^*$ satisfies the flat truncation condition \cite{Nie2013FlatTruncation}, 
which means
\begin{align} \label{cndt:psd:ft}
  \rank\big( M_s(\bar{y}^*) \big) = \rank\big( M_{s-1}(\bar{y}^*) \big)
\end{align}
holds with an integer $1 \leq s \leq k$.
Then we can obtain at least $\rank(M_s(\bar{y}^*))$ minimizers 
of the problem \eqref{opt:CARE:PSD} via the extraction algorithm 
in \cite{Henrion2005extractSolution}. The following algorithm 
is given for solving the problem \eqref{opt:CARE:PSD}.
\begin{alg}
\label{algo:CARE:PSD}

Let $\bar{f}$, $\bar{g}$ and $P$ be as in \eqref{opt:CARE:PSD}.
Set $k := 1$, do the following:

\begin{itemize}

\item[1:]
Solve the SDP problem \eqref{opt:CARE:PSD:Relax}.

\item[2:]
If the relaxation \eqref{opt:CARE:PSD:Relax} is infeasible,
the CARE \eqref{CARE:general} has no positive semi-definite solution and stop.
Otherwise, solve it for a minimizer $\bar{y}^*$. Let $s := 1$.

\item[3:]
Check whether or not $\bar{y}^*$ satisfies the rank condition \eqref{cndt:psd:ft}.
If it is satisfied, extract $\rank(M_s(\bar{y}^*))$ minimizers 
for \eqref{opt:CARE:PSD} and stop.
Otherwise, go to Step 4.

\item[4:]
If $s < k$, let $s := s + 1$ and go to Step 3.
Otherwise, let $k := k + 1$ and go to Step 1.
\end{itemize}
\end{alg}
The MATLAB software packages {\tt YALMIP} \cite{Lofberg2004YALMIP} 
implemented Algorithm \ref{algo:CARE:PSD}.

The dual problem of the problem \eqref{opt:CARE:PSD:Relax} is 
\begin{align} \label{opt:CARE:PSD:Relax:Dual}
  \begin{cases}
    \max & \bar{\gamma}, \\
    \st  & \bar{f}-\bar{\gamma} \in \IQ(\bar{g}, P)_{2k},
  \end{cases}
\end{align}
which is also an SDP problem. Let $\bar{f}_k^*$ be the optimal value 
of the problem \eqref{opt:CARE:PSD:Relax:Dual}. 
In term of weak duality, we have $\bar{f}_k^* \leq \bar{f}_k$ 
for each relaxation order $k$.
Moreover, we have the monotonically increasing sequence 
\[
  \bar{f}_k^*     ~\leq~ \bar{f}_{k+1}^* ~\leq~ 
  \bar{f}_{k+2}^* ~\leq~ \cdots ~\leq~ \bar{f}^*,   
\]
as the relaxation order $k$ increases. 
When the problem \eqref{opt:CARE:PSD:Relax} satisfies the Slater condition, 
i.e., there is an interior point in the feasible set of the problem 
\eqref{opt:CARE:PSD:Relax}, then the problem \eqref{opt:CARE:PSD:Relax:Dual} 
has optimal value and $\bar{f}_k=\bar{f}_k^*$ for each relaxation order $k$.

Under some conditions, there are different upper bounds for 
positive (semi)definite solutions of the CARE \eqref{CARE:general}, i.e., 
we can find a matrix $P_u$ such that $P_u-P_*$ is positive semi-definite 
for the positive semi-definite solution $P_*$ of the CARE \eqref{CARE:general}.

If the CARE \eqref{CARE:general} has an upper bound $P_u$, then we have
\begin{align} \label{care:upperbound}
  \| P_u \|_{F}^2 - \|P\|_{F}^2 \geq 0.
\end{align}
If the constraint \eqref{care:upperbound} is added to 
the problem \eqref{opt:CARE:PSD}, 
then the archimedean condition always holds.

For the problems \eqref{opt:CARE:PSD} and \eqref{opt:CARE:PSD:Relax}, 
the following conclusions are given.
\begin{thm}
Let $\bar{f}^*$ and $\bar{f}_k$ be the optimal value of the problem
\eqref{opt:CARE:PSD} and the relaxation \eqref{opt:CARE:PSD:Relax}, respectively.
Assume that the archimedean condition holds, then

\begin{enumerate}
\item [\rm i)]
The relaxation \eqref{opt:CARE:PSD:Relax} is infeasible for some $k$
if and only if the CARE \eqref{CARE:general} has no positive semi-definite solution.

\item [\rm ii)]
If the relaxation \eqref{opt:CARE:PSD:Relax} has a solution $\bar{y}^*$
and the rank condition \eqref{cndt:psd:ft} holds, then $\bar{f}_k=\bar{f}^*$
and there exists a positive semi-definite solution for the problem \eqref{CARE:general}.
\end{enumerate}
\end{thm}

\begin{pf}
i) Since the problem \eqref{opt:CARE:PSD:Relax} is the relaxation of
the problem \eqref{opt:CARE:PSD}, the relaxation \eqref{opt:CARE:PSD:Relax}
is infeasible, then the problem \eqref{opt:CARE:PSD} is infeasible, i.e.,
the CARE \eqref{CARE:general} has no positive semi-definite solution.

Conversely, if the CARE \eqref{CARE:general} has no positive semi-definite solution,
then the set $\{P\in\sR^{n\times n}~|~\bar{g} = 0, ~P \succeq 0 \}$ is empty. 
Hence, the problem \eqref{opt:CARE:PSD} is infeasible. 
Since the archimedean condition holds,
we have $-1 \in \IQ(\bar{g},P)$. So, $-1 \in \IQ(\bar{g},P)_{2k}$
for all $k$ big enough, which means that \eqref{opt:CARE:PSD:Relax:Dual}
is unbounded from above for all big $k$. By weak duality,
the relaxation \eqref{opt:CARE:PSD:Relax} is infeasible.

ii) The conclusion is from the classical achievements in
\cite{Nie2013FlatTruncation} and \cite{Curto2005TKM}.
\end{pf}

\section {Numerical examples}

In this section, we give several numerical examples to show that 
our approaches are effective. All experiments are performed in 
MATLAB R2020b on a Lenovo laptop with an octa-core CPU @3.20 GHz and RAM 16 GB. 
We use SDP solver {\tt SeDuMi} \cite{Sturm1999SeDuMi} to compute semi-definite 
relaxation problems in Algorithm \ref{algo:CARE:CD} and Algorithm \ref{algo:CARE:PSD}.
The integer $n$ denotes the dimension of 
the coefficient matrix $A$ in the CARE \eqref{CARE:general}.

\begin{exmp}
When $n=3$, consider the CARE \eqref{CARE:general} with
\begin{gather*}
A = \begin{pmatrix}
   1.8 &  0.6 & -0.2 \\ 
   0.8 &  1.6 & -0.2 \\ 
  -0.4 & -0.8 &  2.6 
\end{pmatrix}, \quad
R_1 = \begin{pmatrix} 1 \\ 1 \\ 0 \end{pmatrix}, \\
C = \begin{pmatrix} 1 & 1 & 2 \end{pmatrix}, \quad
R = R_1R_1^{\T}, \quad 
Q = C^{\T}C.
\end{gather*}
Applying Algorithm \ref{algo:CARE:CD} and Algorithm \ref{algo:CARE:PSD},
it is infeasible, which means there is no positive semi-definite solution.
The Newton method, structure-preserving doubling algorithm and 
matrix sign function method are not convergent after 100000 iterations. 
The Schur method is also not convergent.
Therefore, our algorithms are effective for detecting 
the nonexistence of the positive semi-definite solution.
\end{exmp}

\begin{exmp} \label{ex:no-psd-solution:gcare:zwj}
When $n=3$, consider the CARE \eqref{CARE:general} with
\begin{gather*}
A = \begin{bmatrix} 
  1.8 & 0.6 & -0.2 \\ 0.8 & 1.6 & -0.2 \\ -0.4 & -0.8 & 2.6 
\end{bmatrix}, \quad
R = \begin{bmatrix} 1 & 1 & 0 \\ 1 & 1 & 0 \\ 0 & 0 & 0 \end{bmatrix}, \\
Q = \begin{bmatrix} 
  -7.6725 &  6.0546 & -3.8884 \\ 
   6.0546 & -9.6761 &  4.2537 \\
  -3.8884 &  4.2537 & -1.9167 
\end{bmatrix}.
\end{gather*}
Applying Algorithm \ref{algo:CARE:CD} to solve problem \eqref{opt:CARE:CD} 
and Algorithm \ref{algo:CARE:PSD} to solve problem \eqref{opt:CARE:PSD}, 
we obtain the following exact positive semi-definite solutions: 
\begin{gather*}
  P_{*,1} = \begin{bmatrix}
   4.0478 & -2.9288 &  1.2014 \\
  -2.9288 &  3.9168 & -1.0209 \\
   1.2014 & -1.0209 &  0.3888
  \end{bmatrix}, \\ 
  P_{*,2} = \begin{bmatrix}
    4.0886 & -2.8741 &  1.2161 \\
   -2.8741 &  3.9904 & -1.0012 \\
    1.2161 & -1.0012 &  0.3940
  \end{bmatrix},
\end{gather*}
respectively.
After 100000 iterations, the Newton method, 
structure-preserving doubling algorithm 
and matrix sign function method are not convergent. 
The Schur method is also not convergent. Therefore, 
these methods can not find a positive semi-definite solution 
and they are not suitable for solving the CARE.
\end{exmp}

\begin{exmp} \label{ex:no-psd-solution:gcare:zj}
When $n=4$, consider the CARE \eqref{CARE:general} with
\begin{gather*}
A = \begin{pmatrix}
  -1 &  0 & 0 & 0 \\
  -2 &  1 & 0 & 0 \\
   0 &  5 & 5 & 0 \\
   3 & -3 & 3 & 6
\end{pmatrix}, \quad
R = \begin{pmatrix}
   4 & -3 &  0 &  2 \\
  -3 &  7 &  1 &  2 \\
   0 &  1 &  8 & -6 \\
   2 &  2 & -6 &  5
\end{pmatrix}, \\
Q = \begin{pmatrix}
  -3 & -8  & 10 &  15 \\
  -8 & -32 & 4  & -1  \\
  10 &  4  & 20 &  23 \\
  15 & -1  & 23 &  31
\end{pmatrix}.
\end{gather*}
Applying Algorithm \ref{algo:CARE:CD} and Algorithm \ref{algo:CARE:PSD},
it is infeasible, which means there is no positive semi-definite solution.
the Newton method, structure-preserving doubling algorithm, and 
matrix sign function method are not convergent after 100000 iterations. 
The Schur method is also not convergent. 
Therefore, our algorithms are effective for detecting 
the nonexistence of the positive semi-definite solution.
\end{exmp}

\begin{exmp}
(\cite[Example 1]{Liu2020UpperBound})
When $n=3$, consider the CARE \eqref{CARE:general} with
\begin{gather*}
A   = \begin{pmatrix} -1 & 0 & 2 \\ 0 &       -1 & 0 \\ 0 & 0 & -3 \end{pmatrix}, \quad
R_1 = \begin{pmatrix}  1 & 0 & 1 \\ 0 &        1 & 0 \\ 1 & 0 & -1 \end{pmatrix}, \\
C   = \begin{pmatrix}  1 & 0 & 0 \\ 0 & \sqrt{2} & 0 \\ 0 & 0 &  1 \end{pmatrix}, \quad
R = R_1R_1^{\T}, \quad 
Q = C^{\T}C.
\end{gather*}
According to Theorem 2.1 of \cite{Liu2020UpperBound}, 
we obtain the upper bound 
\[
  P_u = \begin{pmatrix}
  0.4730 & 0      & 0 \\
  0      & 0.9459 & 0 \\
  0      & 0      & 0.4730
 \end{pmatrix}.
\]
After adding their corresponding upper bound constraints 
$\| P_u \|_{F}^2 - \|XX^{\T}\|_{F}^2 \geq 0$ and 
$\| P_u \|_{F}^2 - \|P\|_{F}^2 \geq 0$ to problems 
\eqref{opt:CARE:CD} and \eqref{opt:CARE:PSD} respectively,
applying Algorithm \ref{algo:CARE:CD} and Algorithm \ref{algo:CARE:PSD}, 
we obtain the following exact positive semi-definite solution
\[
  P^* = \begin{pmatrix}
  0.3551 & 0.0000 & 0.1372 \\
  0.0000 & 0.7321 & 0.0000 \\
  0.1372 & 0.0000 & 0.2336
 \end{pmatrix}.
\]
All real parts of the eigenvalues of $A-RP^*$ are negative, 
so $P^*$ is a stabilizing solution.
\end{exmp}

\begin{exmp}
(\cite[Example 1]{Zhang2013UpperBound})
When $n=3$, consider the CARE \eqref{CARE:general} with
\begin{gather*}
A = \begin{pmatrix} -3 & 0 & 1 \\ 1 & 2  & 2 \\ 0 & 0 & -3 \end{pmatrix}, ~~
R = \begin{pmatrix}  1 & 3 & 0 \\ 3 & 10 & 1 \\ 0 & 1 &  1 \end{pmatrix}, ~~
Q = \begin{pmatrix}  5 & 0 & 1 \\ 0 & 6  & 0 \\ 1 & 0 & 11 \end{pmatrix}.
\end{gather*}
According to Corollary 3 in \cite{Zhang2013UpperBound}, we obtain the upper bounds 
\[
 P_u = \begin{pmatrix}
  1.1125 & -0.1244 &  0.1635 \\
 -0.1244 &  1.3426 & -0.0476 \\
  0.1635 & -0.0476 & 1.6322
 \end{pmatrix}.
\]
After adding their corresponding upper bound constraints 
$\| P_u \|_{F}^2 - \|XX^{\T}\|_{F}^2 \geq 0$ and 
$\| P_u \|_{F}^2 - \|P\|_{F}^2 \geq 0$ to problems 
\eqref{opt:CARE:CD} and \eqref{opt:CARE:PSD} respectively,
applying Algorithm \ref{algo:CARE:CD} and Algorithm \ref{algo:CARE:PSD}, 
we obtain the following exact positive semi-definite solution
\[
  P^* = \begin{pmatrix}
   0.7655 & -0.1338 &  0.2163 \\
  -0.1338 &  1.0528 & -0.0220 \\
   0.2163 & -0.0220 &  1.5153
 \end{pmatrix}.
\]
All real parts of the eigenvalues of $A-RP^*$ are negative, 
so $P^*$ is a stabilizing solution.
\end{exmp}

\section{Conclusion}

In this paper, we study how to obtain a positive semi-definite solution of the CARE. 
The CARE is formulated as two constrained polynomial optimization problems, 
and we use Lasserre's hierarchy of semi-definite relaxations to solve
polynomial optimization problems. Our approaches can check whether the CARE 
has a positive semi-definite solution or not without strong assumptions. 
If there is a positive semi-definite solution for the CARE, we can obtain 
the exact solution by solving polynomial optimization problems.
Furthermore, if the CARE has an upper bound under some conditions,
our approaches are convergent without the assumption of archimedean condition.
Numerical experiments have shown that our methods are effective.

Lasserre's hierarchy of semi-definite relaxations can check global optimality 
of polynomial optimization problems, but it is computationally expensive.
A natural but challenging work is whether we can find an approach to solve 
the CARE with special structures of coefficient matrices, such as sparse matrices.

\begin{ack}
The work was supported in part by National Key Research 
and Development Plan (2023YFB3001604).
\end{ack}

\end{document}